\documentclass[11pt,a4paper,oneside,english]{amsart}
\usepackage{graphicx}
\usepackage{graphics}
\usepackage{epstopdf}
\usepackage{amsmath} 
\usepackage{amssymb} 
\usepackage{color} 
\usepackage[T1]{fontenc}
\usepackage[latin9]{inputenc}

\usepackage{amsthm}
\usepackage{graphicx}
\usepackage{amssymb}
\usepackage{color}

\usepackage{babel}

\newtheorem{remark}{Remark}

\begin{document}
\title{The stability, persistence and extinction in a stochastic model of the population growth}

\author[A. Korobeinikov  and L. Shaikhet]
{Andrei Korobeinikov and Leonid Shaikhet}

\maketitle


\medskip
{\footnotesize
\centerline{School of Mathematics and Information Science,}
\centerline{Shaanxi Normal University, Xi'an, China}
\centerline{akorobeinikov777@gmail.com}}

\medskip
{\footnotesize
\centerline{Department of Mathematics, Ariel University, Ariel 40700, Israel}
\centerline{leonid.shaikhet@usa.net}}

\medskip

{\centering \textbf{Proposed running head:} Properties of a stochastic population model
\vspace{1cm}\par}

{\centering \textbf{AMS Classification (MSC2010)}\\
 92D30 (primary), 34D20, 60H10 (secondary)\par}

{\centering \vspace{1cm}\par}


\begin{abstract}
In this paper we consider the global qualitative properties of a stochastically perturbed logistic model of population growth. 
In this model, the stochastic perturbations are assumed to be of the white noise type and are proportional to the current population size. 
Using the direct Lyapunov method, we established the global properties of this stochastic differential equation. 
In particular, we found that solutions of the equation oscillate around an interval, and explicitly found the end points of this interval. 
Moreover, we found that, if the magnitude of the noise exceeds a certain critical level (which is also explicitly found), then 
the stochastic stabilisation (``stabilisation by noise'') of the zero solution occurs. In this case,  (i) the origin is the lower boundary 
of the interval, and (ii) the extinction of the population due to stochasticity occurs almost sure (a.s.) for a finite time. 

\medskip
{\bf Keywords:} logistic differential equation, stability, extinction, persistence of a population, stochastic perturbations, Wiener process, 
Lyapunov function, the direct Lyapunov method.
\end{abstract}

\pagebreak


The growth of a population in an environment with a limited currying capacity is usually described by the logistic differential equation
\begin{equation}\begin{array}{l}
\label{eq1}
\dot x(t)=ax(t)-bx^2(t),
\end{array}\end{equation}
where $a$ and $b$ are positive parameters, and ratio $K=\dfrac{a}{b}$ is the currying capacity of the environment~\cite{Murrey}.
This equation is defined for all $x\ge0$ and has two equilibria: an unstable equilibrium state at the origin
and asymptotically stable positive equilibrium state $x^*=\dfrac{a}{b}$.
The stability of equilibrium state $x^*$ is easy to verify: indeed, substituting $x(t)=y(t)+x^*$ into \eqref{eq1},
for $y(t)$ we obtain equation $\dot y(t)=-ay(t)-by^2(t)$.

Let us suppose now that  equation \eqref{eq1} is exposed to stochastic perturbations that are of the type
of white noise and are proportional to the current population size. Then the ordinary differential  equation \eqref{eq1} transforms
to the following Ito's stochastic differential equation \cite{GiSk}:
\begin{equation}\begin{array}{l}\label{eq1s}
dx(t)=(ax(t)-bx^2(t))dt+\sigma x(t)dw(t),
\end{array}\end{equation}
where $\sigma$ is a constant and $w(t)$ is the standard Wiener process. Please note that the stochastic perturbations
proportional to the current population size appear to be the most natural in biological systems. However,
the analysis of such equations can be challenging as this kind of equations has no positive equilibria.
Indeed, it is easy to see that the origin  remains to be an unstable equilibrium state of the model.
At the same time, for the stochastic equation \eqref{eq1s} the point $x^*$ is not an equilibrium state anymore,
and even is not a solution, hence, investigating its properties may appear to be meaningless.
Nonetheless, the direct Lyapunov method enables us to investigate the behaviour of solutions of stochastically
perturbed logistic equation \eqref{eq1s} in the neighbourhood of the point $x^*$
(the positive equilibrium state of the deterministic model \eqref{eq1}).

We consider the Lyapunov function
\begin{equation}\begin{array}{l}\label{v}
v(x)=\dfrac{x}{x^*}-\ln\dfrac{x}{x^*}-1,
\end{array}\end{equation}
which is defined for all $ x>0 $.
This function was invented by Volterra \cite{Volt} and is extremely useful in mathematical biology.

Let $L$ be the generator of equation \eqref{eq1s} \cite{GiSk}. Using $x^*=\dfrac{a}{b}$, we have
\begin{equation}\begin{array}{l}\label{Lv1}\aligned
Lv(x)=&\left(\dfrac{1}{x^*}-\dfrac{1}{x}\right)(ax(t)-bx^2(t))+\dfrac{1}{2}\sigma^2\\
=& b \left( \dfrac{x(t)}{x^*}-1 \right)(x^* - x(t))+\dfrac{1}{2}\sigma^2\\
=& - \dfrac{b}{x^*} (x(t)-x^*)^2+\dfrac{1}{2}\sigma^2.
\endaligned\end{array}\end{equation}
Hence, $Lv(x*)=\dfrac{1}{2}\sigma^2$ and $Lv(x)=0$ at points $x_1$ and $x_2$, where equality
$$
(x(t)-x^*)^2 = \dfrac{x^*}{2b} \sigma^2
$$
holds. It is easy to see that, if $\sigma^2<2a$, then $x_1=x^*\left(1-\dfrac{|\sigma|}{\sqrt{2a}}\right)$.
If $\sigma ^2 \geq 2a$, then $x_1$ is undefined (out of the region $(0, +\infty)$, where the Lyapunov function is defined).
For convenience, in this case we set  $x_1=0$. Point $x_2=x^*\left(1+\dfrac{|\sigma|}{\sqrt{2a}}\right)$
exists for all $\sigma$. $Lv(x)\in\left(0,\dfrac{1}{2}\sigma^2\right]$ for $x\in X=(x_1,x_2)$, whereas
$Lv(x) < 0$  for all $x \in (0,x_1)$ and $x \in (x_2,\infty)$.

Note that the Lyapunov function \eqref{v} decreases for $x\in(0,x^*)$ and increases for all $x>x^*$. Besides,~\cite{GiSk}
\begin{equation}\begin{array}{l}\label{ELv1}\aligned
\mathbf Ev(x(t))-\mathbf Ev(x(0))=\int^t_0\mathbf ELv(x(s))ds,
\endaligned\end{array}\end{equation}
and, hence, if $Lv(x(s))>0$ for $s\in(0,t)$ then $\mathbf Ev(x(t))>\mathbf Ev(x(0))$.
This fact allows us to make decisive conclusions regarding  system long-term behaviour.
Indeed, this fact implies that, if $x(0)\in(x_1,x^*)$, then $x(t)$ decreases almost sure (a.s.) to $x_1$,
and, if $x(0)\in(x^*,x_2)$, then $x(t)$ increases a.s. to $x_2$.
Likewise, if $x(0)<x_1$ then $Lv(x(s))<0$ for $s\in(0,t)$ and $\mathbf Ev(x(t))<\mathbf Ev(x(0))$.
It means that $x(t)$ increases a.s. to $x_1$.
Finally, if $x(0)>x_2$ then $Lv(x(s))<0$ for $s\in(0,t)$ and $\mathbf Ev(x(t))<\mathbf Ev(x(0))$.
It means that $x(t)$ decreases a.s. to $x_2$.
This implies that for $\sigma^2<2a$ and all positive initial conditions the solutions oscillate around the interval $(x_1,x_2)$.

Another important conclusion can be made regarding a possibility of the stochastic extinction of the population:
it immediately follows from our analysis that if $\sigma^2\ge2a$ then $x_1=0$, and, hence, the extinction
due to stochastic perturbations is certain (a.s.) and occurs for a finite time. 
This conclusion implies that for $\sigma^2 \geq 2a$  the origin should be a stable equilibrium state and, hence, 
the so-called ``stochastic stabilisation'', or ``stabilisation by noise'', of the zero solution occurs at $\sigma^2 = 2a$. 
Indeed, it is easy to see that the stochastic stabilisation occurs.
For equation \eqref{eq1s}, Lyapunov function $V(x)=|x|^{1-\frac{2a}{\sigma^2}}$ with condition $2a<\sigma^2$ 
from~\cite{Khas} yields
\begin{equation}\begin{array}{l}\label{Lvk}\aligned
LV(x)=&\left(1-\dfrac{2a}{\sigma^2}\right)|x|^{-\frac{2a}{\sigma^2}}(ax-bx^2)
-\frac{1}{2}\sigma^2x^2\frac{2a}{\sigma^2}\left(1-\frac{2a}{\sigma^2}\right)|x|^{-\frac{2a}{\sigma^2}-1}\\
=&-\left(1-\frac{2a}{\sigma^2}\right)b|x|^{2-\frac{2a}{\sigma^2}}\le0.
\endaligned\end{array}\end{equation}
Therefore, for arbitrary $b\ge0$ and for all  $\sigma^2>2a $, the zero solution of equation \eqref{eq1s} is stable in probability. 
In the case $b=0$ it is the classical Khasminskii's example of the stochastic stabilisation~\cite{Khas}. 
\begin{remark} 
The instability of the zero solution for all $\sigma^2 < 2a$ follows from analysis of the same Lyapunov function. 
\end{remark}

\begin{remark}
Please note, that if $\sigma=0$, then $x_1=x_2=x^*$, and $x^*$ is globally asymptotically stable equilibrium state.
\end{remark}

Figures 1 -- 5 illustrate these analytic conclusions.
Fig.~1 shows trajectories of solutions of the equation \eqref{eq1s} with $a=1.5$, $b=1$ and $\sigma=0.25$
for $x(0)=2.3$ (10 blue lines) and for $x(0)=0.65$ (10 green lines). It is easy to see that for each of these
initial conditions all trajectories are located essentially between the lines $x_1<x^*=1.5<x_2$ (represented by red dashed lines).
Fig.~2 shows trajectories for $\sigma=0.025$ and for the same values of all other parameters and the same initial conditions.
Fig.~3 shows a trajectory with $x(0)=2$ (the blue line) and a trajectory with $x(0)=0.1$ (the green line) for $a=2.5$, $b=1$ and $\sigma=1.5$.
Fig.~4 depicts the solutions with initial conditions taken precisely at the points $x_1=1.07$ (the green line) and $x_2=1.93$ (the blue line)
for $a=1.5$, $b=1$ and $\sigma=0.5$. It is easy to see that the results shown in these figures coincide with our analytical conclusions.
Fig.~5 shows an example of the stochastic extinction. In this figure, a solution of equation \eqref{eq1s}  with 
$a=b=1$, $\sigma=2.45$, $x_1=0$, $x_2=2.73$, $x(0)=1.75$ is shown.

\begin{figure}
\begin{center}
\includegraphics[height=5cm]{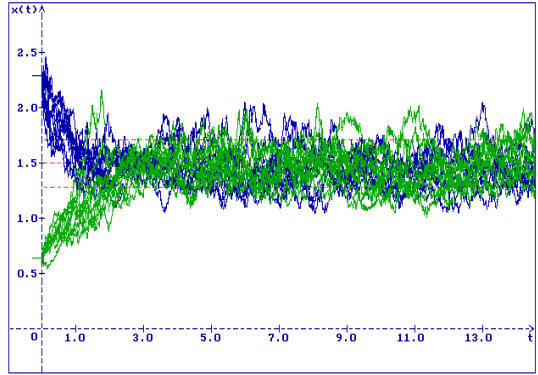}
\caption{Trajectories of solutions of the equation \eqref{eq1s} with $a=1.5$, $b=1$ and $\sigma=0.25$
for $x(0)=2.3$ (10 blue lines) and for $x(0)=0.65$ (10 green lines).
Red dashed lines corresponds to $x_1=1.28$, $x^*=1.5$ and $x_2=1.72$.}
\label{Fig1}
\end{center}
\end{figure}

\begin{figure}
\begin{center}
\includegraphics[height=5cm]{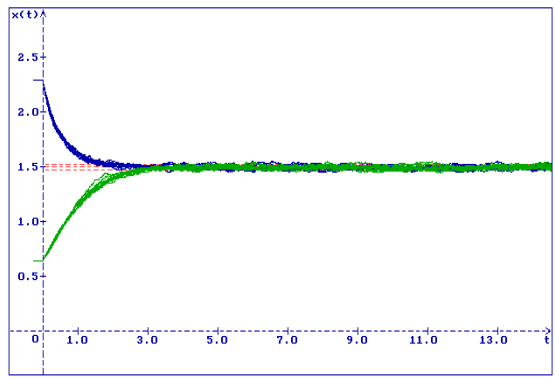}
\caption{Trajectories for $\sigma=0.025$, $x_1=1.47$, $x_2=1.52$; the other parameters are the same as in Fig.~1. }
\label{Fig2}
\end{center}
\end{figure}

\begin{figure}
\begin{center}
\includegraphics[height=5cm]{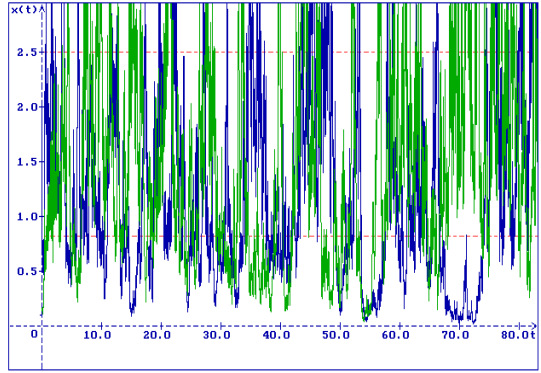}
\caption{Trajectories of solutions of the equation \eqref{eq1s} with $a=2.5$, $b=1$ and $\sigma=1.5$
for $x(0)=2$ (blue line) and for $x(0)=0.1$ (green line). Red dashed lines corresponds to $x_1=0.82$, $x^*=2.5$ and $x_2=4.18$. }
\label{Fig3}
\end{center}
\end{figure}

\begin{figure}
\begin{center}
\includegraphics[height=5cm]{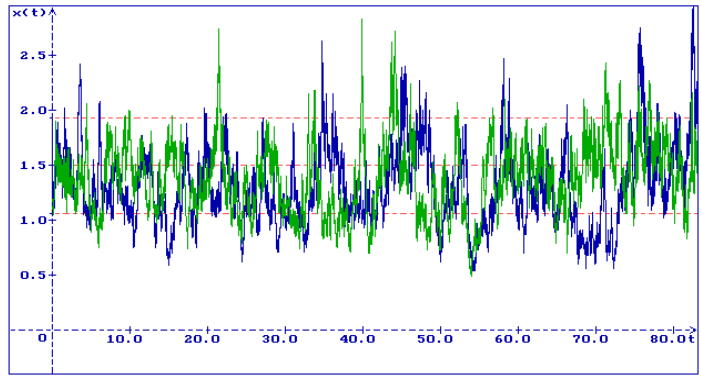}
\caption{Trajectories of solutions of the equation \eqref{eq1s} with $a=1.5$, $b=1$ and $\sigma=0.5$, with initial conditions taken precisely at the points
$x_1=1.07$ (green) and $x_2=1.93$ (blue) . }
\label{Fig4}
\end{center}
\end{figure}


\begin{figure}
\begin{center}
\includegraphics[height=5cm]{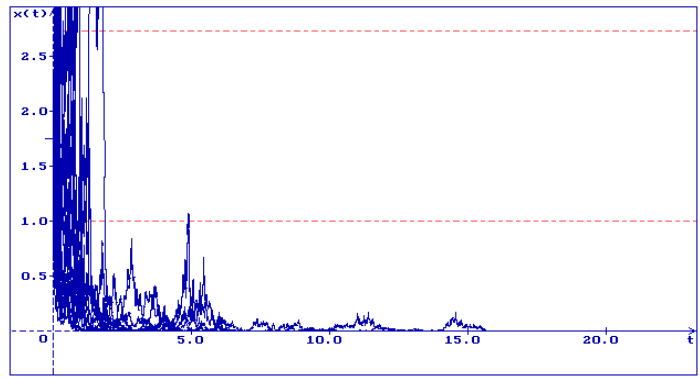}
\caption{A case of the stochastic extinction: $a=b=1$, $\sigma=2.45$, $x_1=0$, $x_2=2.73$, $x(0)=1.75$}
\label{Fig6}
\end{center}
\end{figure}

\end{document}